\documentclass[12pt]{amsart}

\usepackage{mathrsfs}

\usepackage{amsmath,amsthm,amssymb}
\font\teneufm=eufm10
\font\seveneufm=eufm7
\font\fiveeufm=eufm5
\newfam\frakturfam

\textfont\frakturfam=\teneufm
\scriptfont\frakturfam=\seveneufm
\scriptscriptfont\frakturfam=\fiveeufm

\newtheorem{pr}{Proposition}

\newtheorem{lm}{Lemma}

\newtheorem{co}{Corollary}

\newtheorem{prob}{Problem}
\def\bee{\begin{eqnarray}}
\def\bes{\begin{eqnarray*}}
\def\eee{\end{eqnarray}}
\def\ees{\end{eqnarray*}}

\def\a{\alpha}

\def\Proof{{\sl Proof.}\ }

\title{Left-symmetric algebras of derivations of free algebras}

\begin{document}
\date{}
\maketitle

\begin{center}

{\bf Ualbai Umirbaev}\footnote{Supported by a grant CNPq 401287/2012-2; Eurasian National University,
 Astana, Kazakhstan and
 Wayne State University,
Detroit, MI 48202, USA,
e-mail: {\em umirbaev@math.wayne.edu}}

\end{center}

\begin{abstract} A structure of a left-symmetric algebra on the set of all derivations of  a free algebra is introduced such that its commutator algebra becomes the usual Lie algebra of derivations. Left and right nilpotent elements of left-symmetric algebras of derivations are studied. Simple left-symmetric algebras of derivations and Novikov algebras of derivations are described. It is also proved that the positive part of the left-symmetric algebra of derivations of a free nonassociative symmetric $m$-ary algebra in one free variable is generated by one derivation and some right nilpotent derivations are described.
\end{abstract}


\noindent


\section{Introduction}

\hspace*{\parindent}

If $A$ is an arbitrary algebra over a field $k$, then the set $\mathrm{Der}_kA$ of all $k$-linear derivations  of $A$ forms a Lie algebra. If $A$ is a free algebra, then it is possible to define a multiplication  $\cdot$ on $\mathrm{Der}_kA$ such that it becomes a left-symmetric algebra and its commutator algebra becomes the Lie algebra $\mathrm{Der}_kA$ of all derivations of $A$. The language of the left-symmetric algebras of derivations is more convenient to describe some combinatorial properties of derivations.

Recall that an algebra $B$  over $k$ is called {\em left-symmetric} \cite{Burde2006} if $B$ satisfies the identity
\bee\label{f1}
(xy)z-x(yz)=(yx)z-y(xz).
\eee
This means that the associator $(x,y,z):=(xy)z-x(yz)$ is symmetric with respect to two left arguments, i.e.,
\bes
(x,y,z)=(y,x,z).
\ees
The variety of
left-symmetric algebras is Lie-admissible, i.e., each
left-symmetric algebra $B$ with the operation $[x,y]:=xy-yx$ is a
Lie algebra. A linear basis of free left-symmetric algebras is constructed in \cite{Segal}.

In this paper we introduce the left-symmetric algebras of derivations $\mathscr{L}(A)$ of an arbitrary free algebra $A$ and prove some introductory results on their structure. Simple left-symplectic algebras of derivations are described.

We prove that the set of all left nilpotent elements of $\mathscr{L}(A)$ coincides with the set  of all locally nilpotent derivations of $A$. Thus, a well known combinatorial question on the structure of locally nilpotent derivations of free algebras becomes a purely structural algebraic question of left-symmetric algebras of derivations. We also show that if the Jacobian matrix $J(D)$ of a derivation $D\in \mathscr{L}(A)$ is nilpotent then the right multiplication operator $R_D$ is nilpotent and show that the converse is true in many cases, in particular, for free algebras of the Nielson-Schreier varieties of algebras.

The remaining part of the paper is devoted to examples of free one generated algebras $A$ such that the positive part $\mathscr{L}(A)_+$ of $\mathscr{L}(A)$ is generated by one element. In particular, Novikov algebras of derivations are described.  We also prove that the positive part of the left-symmetric algebra of derivations of a free nonassociative symmetric $m$-ary one generated algebra is generated by one element. We prove that some right nilpotent derivations have nilpotent Jacobian matrices.
The study of $m$-ary symmetric algebras is motivated by an approach to the Jacobian Conjecture by L.M. Dru\.zkowski and K. Rusek \cite{DR85}, G. Gorni and G. Zampieri \cite{GZ96}, and A.V. Yagzhev \cite{Yagzhev00-1} (see also \cite{BBRY,UU2014-2}), which reduces the study of the Jacobian Conjecture to the study of symmetric ternary algebras with an Engel identity. In fact, the left-symmetric algebras of derivations are introduced as a convenient language of the study of $m$-ary symmetric algebras.

This paper is organized as follows. In Section 2 we define left-symmetric algebras of derivations and describe simple algebras of derivations. Section 3 is devoted to the study of left and right nilpotent elements and presents generalizations of some results from \cite{UU2014-3}. Section 4 is devoted to Novikov algebras of derivations and Section 5 is devoted to $m$-ary symmetric algebras.

\section{Left-symmetric product}

\hspace*{\parindent}

Taking  some applications into account, we consider algebras with one $m$-ary  multilinear operation $\langle \cdot,\cdot,\ldots,\cdot\rangle$ ($m\geq 2$) over a field $k$.
So, let ${\mathfrak M}$ be an arbitrary variety of $m$-ary algebras over $k$ and let $A=k_{\mathfrak M}\langle x_1,x_2,\ldots,x_n\rangle$ be a free algebra of ${\mathfrak M}$ freely generated by $x_1,x_2,\ldots,x_n$. If ${\mathfrak M}$ is a unitary variety of binary algebras, then we assume that $A$ contains an identity element. For any $n$-tuple $F=(f_1,f_2,\ldots,f_n)'$ (column) of elements of $A$, denote by
\bes
D_F=f_1\partial_1+f_2\partial_2+\ldots+f_n\partial_n
\ees
the derivation of $A$ defined by $D_F(x_i)=f_i$ for all $i$. Notice that this form is unique for any derivation $D\in \mathrm{Der}_kA$. For any $u=a\partial_i$ and $v=b\partial_j$, where $a,b\in A$, put
\bes
u\cdot v= ((a\partial_i)(b))\partial_j.
\ees
Extending this operation by distributivity, we get a well defined bilinear operation $\cdot$ on $\mathrm{Der}_kA$.
Denote this algebra by $\mathscr{L}(A)$.
\begin{lm}\label{l1} The following statements are true:

(i) $\mathscr{L}(A)$ is a left-symmetric algebra;

(ii) The commutator algebra of $\mathscr{L}(A)$ is the Lie algebra $\mathrm{Der}_kA$ of all derivations of $A$.
\end{lm}
\Proof Let $x,y\in \mathscr{L}(A)$. Denote by $[x,y]=x\cdot y-y\cdot x$ the commutator in $\mathscr{L}(A)$. For now denote by $\{x,y\}$ their product in $\mathrm{Der}_kA$. In order to prove (ii), we have to check that
\bes
[x,y](a)=\{x,y\}(a)
\ees
 for all $a\in A$. Notice that
\bes
\{x,y\}(a)=x(y(a))-y(x(a))
\ees
by definition. Taking into account that $[x,y]$ and $\{x,y\}$ are both derivations, we can assume that $a=x_t$.
So, it is sufficient to check that
\bes
(x\cdot y-y\cdot x)(x_t)=x(y(x_t))-y(x(x_t)).
\ees
We may also assume that $x=u\partial_i$ and $y=v\partial_j$. If $t\neq i,j$, then all components of the last equality are zeroes. If $t=i\neq j$ or $t=i=j$, then it is also true. This proves (ii).

Assume that $x,y\in \mathscr{L}(A)$ and $z=a\partial_t$. Then
\bes
(x,y,z)=(xy)z-x(yz)=[(xy)(a)-x(y(a))]\partial_t,\\
(y,x,z)=(yx)z-y(xz)=[(yx)(a)-y(x(a))]\partial_t.
\ees
To prove (\ref{f1}) it is sufficient to check that
\bes
[x,y](a)=x(y(a))-y(x(a))=\{x,y\}(a).
\ees
This follows immediately from (ii). $\Box$

The algebra $\mathscr{L}(A)$ is called the {\em left-symmetric algebra of all derivations} of the free algebra $A$.

Consider the grading
\bes
A=A_0\oplus A_1\oplus A_2\oplus\ldots\oplus A_s\oplus\ldots,
\ees
where $A_i$ is the space of all homogeneous elements of degree $i\geq 0$. Notice that $A_0=k1$ if and only if $m=2$ and ${\mathfrak M}$ is a unitary variety of algebras. Otherwise $A_0=0$. If $m\geq 3$ and $i\neq l(m-1)+1$, then $A_i=0$.

The left-symmetric algebra $\mathscr{L}(A)$ has a natural grading
\bee\label{f2}
\mathscr{L}(A)=L_{-1}\oplus L_{0}\oplus L_{1}\oplus \ldots \oplus L_{s}\oplus\ldots,
\eee
where $L_{i}$ is the linear span of elements of the form $a\partial_j$ with $a\in A_{i+1}$ and $1\leq j\leq n$. Elements of $L_{s}$ are called {\em homogeneous} derivations of $A$ of degree $s$.

If $A$ has an identity element, then $L_{-1}=k\partial_1+\ldots+k\partial_n$. Otherwise $L_{-1}=0$. The space $L_{0}$ is a subalgebra of $\mathscr{L}(A)$ and is isomorphic to the matrix algebra $M_n(k)$. The element
\bes
D_X=x_1\partial_1+x_2\partial_2+\ldots+x_n\partial_n
\ees
is the identity element of $L_{0}$ and is the right identity element of $\mathscr{L}(A)$. The left-symmetric algebra
$\mathscr{L}(A)$ has no identity element.

The subalgebra
\bes
\mathscr{L}(A)_+=L_1\oplus L_2\oplus \ldots \oplus L_n\oplus\ldots
\ees
of $\mathscr{L}(A)$ will be called {\em the positive part} of $\mathscr{L}(A)$ with respect to the grading (\ref{f2}).

{\bf Example 1.} The best known example of a left-symmetric algebra of derivations is the
 left-symmetric Witt algebra
 $\mathscr{L}_n$ of index $n$ \cite{Burde2006,UU2014-3}, i.e., the left-symmetric algebra of all  derivations of the polynomial algebra $P_n=k[x_1,x_2,\ldots,x_n]$.
 The commutator algebra of $\mathscr{L}_n$ is the Witt algebra $W_n$ of index $n$.

Consider the grading
\bes
P_n=T_0\oplus T_1\oplus T_2\oplus\ldots T_s\oplus\ldots
\ees
 of the polynomial algebra $P_n=k[x_1,x_2,\ldots,x_n]$, where $T_i$ is the space of all homogeneous polynomials of degree $i$. Then
\bes
\mathscr{L}_n=Q_{-1}\oplus Q_{0}\oplus Q_{1}\oplus \ldots \oplus Q_{s}\oplus\ldots,
\ees
where $Q_{i}$ is the linear span of all elements of the form $a\partial_j$ with $a\in T_{i+1}$ and $1\leq j\leq n$.

{\bf Example 2.} If $m\geq 3$, then the subalgebra
\bes
\mathscr{L}_{n,m}=Q_{0}\oplus Q_{m-1}\oplus \ldots \oplus Q_{ms-s}\oplus\ldots
\ees
of $\mathscr{L}_n$ is a left-symmetric algebra of derivations of a free algebra.

Define an $m$-ary operation on the polynomial algebra $P_n$ by
\bes
\langle a_1,a_2,\ldots,a_m\rangle= a_1a_2\ldots a_m, \ \ \ a_1,a_2,\ldots,a_m\in P_n.
\ees
Denote by $P_{n,m}$ the subalgebra of  $P_n$ with respect to this $m$-ary operation generated by all $x_i$, where $1\leq i\leq n$. Obviously, $P_{n,2}$ coincides with $P_n$. If $m\geq 3$, then
\bes
P_{n,m}=T_1\oplus T_m\oplus\ldots T_{sm-s+1}\oplus\ldots.
\ees
It is not difficult to show that $P_{n,m}$ is a free $m$-ary algebra freely generated by $x_1,x_2,\ldots,x_n$ and $\mathscr{L}_{n,m}\simeq \mathscr{L}(P_{n,m})$.

{\bf Example 3.} Let $B=kx_1\oplus kx_2\oplus\ldots kx_n$ be an $n$-dimensional algebra with trivial multiplication. Then $\mathscr{L}(B)\simeq M_n(k)$.

We show that $\mathscr{L}_n$ and $M_n(k)$ are the only simple left-symmetric algebras of derivations.

\begin{pr}\label{p1} Let $A$ be an arbitrary $m$-ary free algebra in the variables $x_1,x_2,\ldots,x_n$ over a field $k$ of characteristic zero. Then the left-symmetric algebra $\mathscr{L}(A)$ is simple if and only if $\mathscr{L}(A)$  is isomorphic to the left symmetric algebra $\mathscr{L}_n$ or associative matrix algebra $M_n(k)$.
\end{pr}
\Proof  Suppose that $\mathscr{L}(A)$ is simple. First consider the case when $A$ does not contain an identity element.
Then $L_{-1}=0$ in the grading (\ref{f2}). We have $L_{0}=M_n(k)$. Obviously, $\mathscr{L}(A)$ is simple if and only if $L_i=0$ for all $i\geq 0$. Notice that this implies $\langle A,\ldots,A\rangle=0$, i.e., $A$ is an $n$-dimensional vector space with zero multiplication.

Assume that $A$ is a binary algebra with an identity element. Let $I$ be a nontrivial verbal ideal of $A$, i.e., $0\neq I\neq A$ and $\phi(I)\subseteq I$ for every endomorphism $\phi$ of $A$. Then $I$ is differentially closed and the set of all derivations
of the form
\bes
f_1\partial_1+f_2\partial_2+\ldots+f_n\partial_n, \ \ \ f_i\in I,
\ees
is a nontrivial ideal of $\mathscr{L}(A)$. This implies that $A$ does not have nontrivial verbally closed ideals.

Notice that every $T$-ideal is verbally closed \cite{KBKA}. Let $J$ be the associator ideal of $A$ and $K$ be the commutator ideal of $A$, i.e., $J$ is generated by all associators $(x,y,z)=(xy)z-x(yz)$ and $K$ is generated by all commutators $[x,y]=xy-yx$. Notice that $J$ and $K$ are both verbally closed. We have also $J\neq A$ and $K\neq A$ since $A$ is homogeneous and $1\notin J,K$. Consequently, $J=K=0$ since $A$ has no nontrivial verbally closed ideals.
This means that $A$ is an associative and commutative algebra. If $A$ is a free algebra of the variety of all associative and commutative algebras, then $A\simeq k[x_1,x_2,\ldots,x_n]$.
In this case $\mathscr{L}(A)\simeq \mathscr{L}_n$. Otherwise $A$ belongs to a proper subvariety of all associative and commutative algebras. It is well-known that every proper subvariety of all associative and commutative algebras is nilpotent. But $A$ is not nilpotent since it has an identity element. $\Box$

\section{Nilpotent elements}

\hspace*{\parindent}

Let $G$ be an arbitrary binary (nonassociative) algebra.
  If $a\in G$, then put $a^1=a^{[1]}=a$, $a^{r+1}=a(a^r)$, and $a^{[r+1]}=(a^{[r]})a$ for any $r\geq 1$.
It is natural to say that $a$ is left nilpotent if $a^m=0$ for some $m\geq 2$. Similarly,  $a$ is right nilpotent if $a^{[m]}=0$ for some $m\geq 2$.

Recall that ${\mathfrak M}$ is an arbitrary variety of $m$-ary algebras over $k$ and  $A=k_{\mathfrak M}\langle x_1,x_2,\ldots,x_n\rangle$ is a free algebra of ${\mathfrak M}$.

\begin{lm}\label{l2}
A derivation $D$ of $A$ is locally nilpotent if and only if $D$ is a left nilpotent element of $\mathscr{L}(A)$.
\end{lm}
\Proof It is well known that a derivation $D$ of $A$ is locally nilpotent if and only there exists a positive integer $p$ such that $D$ applied $p$ times to $x_i$ gives zero for all $i$.

Suppose that $D=D_F$ and put
\bes
H_i=\underbrace{D(D\ldots(D(D}_iX))\ldots)
\ees
 for all $i\geq 1$. Note that $H_1=F$. We have also
\bee\label{f3}
D_F D_G= D_{D_F(G)}
\eee
by the definition of the left symmetric product.
 By induction on $i$ and (\ref{f3}), it is easy to show that $D^i=D_{H_i}$. Consequently, $D^m=0$ if and only if $H_m=0$. Note that $H_m=0$ means that $D$ applied $m$ times to $x_i$  gives $0$ for all $i$. $\Box$

 Let $B$ be an arbitrary $m$-ary algebra.
For any $b_1,b_2,\ldots,b_{m-1}\in B$ and $1\leq i\leq m$ denote by
\bes
M_i(b_1,b_2,\ldots,b_{m-1}): B \longrightarrow B
\ees
 the linear operator on $B$ defined by
\bes
M_i(b_1,b_2,\ldots,b_{m-1})(x)=\langle b_1,\ldots,b_{i-1},x,b_{i+1},\ldots,b_{m-1}\rangle
\ees
for all $x\in B$. If $m=2$, then $M_1(b)$ becomes the operator of right multiplication $R_b$ and $M_2(b)$ becomes the operator of left multiplication $L_b$. The associative algebra $M(B)$ (with identity) generated by all $M_i(b_1,b_2,\ldots,b_{m-1})$, where $b_1,b_2,\ldots,b_{m-1}\in B$ and $1\leq i\leq m$, is called the {\em multiplication algebra} of $B$ \cite{KBKA}.

Following \cite{Um11,Um}, we give a short definition of the universal (multiplicative) enveloping algebra $A^e$  and the universal derivation $\Omega: A \rightarrow \Omega_A$ of the free algebra $A$. Consider
$C=k_{\mathfrak M}\langle x_1,x_2,\ldots,x_n,y_1,y_2,\ldots,y_n\rangle$, where $y_1,y_2,\ldots,y_n$ are new variables. Denote by $\Omega_A$ the subspace of all homogeneous elements of $C$ of degree $1$ with respect to the variables  $y_1,y_2,\ldots,y_n$. Denote by $\Omega: C \rightarrow C$ the derivation of $C$ defined by $x_i\mapsto y_i, y_i\mapsto 0$ for all $i$. Notice that $\Omega(A)\subseteq \Omega_A$.
The derivation
\bes
\Omega: A \rightarrow \Omega_A
\ees
 is called the {\em universal} derivation of $A$.

For any $b_1,b_2,\ldots,b_{m-1}\in A$ and $1\leq i\leq m$ we define the operator
\bes
U_i(b_1,b_2,\ldots,b_{m-1}): \Omega_A \longrightarrow \Omega_A
\ees
as the restriction of the multiplication operator $M_i(b_1,b_2,\ldots,b_{m-1})\in M(C)$ to $\Omega_A$.

The {\em universal enveloping} algebra $A^e$ of $A$ is the associative $k$-algebra (with identity) generated by all $U_i(b_1,b_2,\ldots,b_{m-1})$, where $b_1,b_2,\ldots,b_{m-1}\in A$ and $1\leq i\leq m$.

By definition, $\Omega_A$ is a free left $A^e$-module generated by $y_i$, $1\leq i\leq n$. For every $b\in A$ there exist unique elements $u_1,u_2,\ldots,u_n\in A^e$ such that
\bes
\Omega(b)=u_1y_1+u_2y_2+\ldots+u_ny_n.
\ees
The elements $u_i=\frac{\partial b}{\partial x_i}$ are called the Fox derivatives of $b\in A$ \cite{Um11,Um}.

Every $n$-tuple $F=(f_1,f_2,\ldots,f_n)'$ of elements of $A$ represents the endomorphism $F$ of $A$ such that $F(x_i)=f_i$ for all $i$. Denote by $J(F)=(\partial_j(f_i))_{1\leq i,j\leq n}$  the Jacobian matrix of $F$. An analogue of the Jacobian Conjecture for the free algebra $A$ can be formulated as follows: If $J(F)$ is invertible over $A^e$, then $F$ is an automorphism of $A$. This conjecture is true for free associative algebras \cite{Dicks,Schofield}, for free Lie algebras and superalgebras \cite{Reutenauer,Shpilrain,Umi9,ZM}, and for free nonassociative, commutative (characteristic $\neq 2$), and anticommutative algebras \cite{Yagzhev2} (see also \cite{Um11,Um}).

There is a natural homomorphism
\bee\label{f4}
\varphi: A^e\rightarrow M(A)
\eee
defined by $U_i(b_1,b_2,\ldots,b_{m-1})\mapsto M_i(b_1,b_2,\ldots,b_{m-1})$ for all
$b_1,b_2,\ldots,b_{m-1}\in A$ and $1\leq i\leq m$. This homomorphism is not necessarily an isomorphism. For example, if ${\mathfrak L}$ is the variety of all Lie algebras and $L$ is a free algebra in one variable $x$, then $L^e=k[l_x]$ is a polynomial algebra in $l_x$ and $M(L)=k$. This homomorphism turns $A$ into an  $A^e$-module.

Notice that every derivation $D$ of $A$ has the form $D=D_F$ for some $n$-tuple $F$. We put $J(D)=J(F)$. So, the Jacobian matrix of every derivation of $A$ is defined.

\begin{lm}\label{l3} Let $F$ and $G$ be two arbitrary $n$-tuples of elements of $A$.
Then
\bes
D_F D_G=D_{J(G)F}=D_{J(D_G)F}.
\ees
\end{lm}
\Proof Notice that for any $h\in A$ we have
\bes
D_F(h)=\sum_{i=1}^n \frac{\partial h}{\partial x_i}y_i|_{y_i:=f_i}
=(\frac{\partial h}{\partial x_1},\ldots,\frac{\partial h}{\partial x_n})F.
\ees
Using this equality and (\ref{f3}), we get
\bes
D_F D_G=D_{D_F(G)}=D_{J(G)F}.      \  \  \  \ \ \  \ \Box
\ees

\begin{co}\label{c1}
Let $D$ be an arbitrary element of $\mathscr{L}(A)$. If the Jacobian matrix $J(D)$ of $D$ is nilpotent, then $R_D$ is a nilpotent element of $M(\mathscr{L}(A))$.
\end{co}
\Proof It follows from Lemma \ref{l3} that
\bes
(\ldots((D_F\underbrace{D_G)D_G)\ldots )D_G}_m=D_{J(G)^mF}.
\ees
This proves the lemma. $\Box$

\begin{co}\label{c2} If the Jacobian matrix $J(D)$ of $D\in\mathscr{L}(A)$ is nilpotent, then $D$ is a right nilpotent element of $\mathscr{L}(A)$.
\end{co}

Is the converse of the statement of Corollary \ref{c1} true? In other words: Does the equality $J(G)^mF=0$ for all $F$ imply that $J(G)^m=0$? It is true for polynomial algebras \cite{UU2014-3}. In the general case, it is a delicate question related to the study of the homomorphism (\ref{f4}). We prove the converse of Corollary \ref{c1} for some varieties of algebras, including the varieties of Nielson-Schreier algebras \cite{MSY}. Recall that a variety of algebras ${\mathfrak M}$ is called Nielsen-Schreier if every subalgebra of a free algebra is free in ${\mathfrak M}$.

We give some conditions for the homomorphism (\ref{f1}) to be an isomorphism.
The least natural number $n$ such that the variety $Var(k_\mathfrak{M}<x_1,x_2,\ldots,x_n>)$ of algebras generated by  $k_\mathfrak{M}<x_1,x_2,\ldots,x_n>$ is equal to $\mathfrak{M}$ is called the {\em base rank} $\mathrm{rb}(\mathfrak{M})$ of the variety $\mathfrak{M}$ \cite{Shestakov77}. The base rank of the variety of associative \cite{Malcev52} and Lie algebras \cite{Shirshov58} is $2$, Novikov \cite{MLU2011-1} and dual Leibniz \cite{NU} algebras is $1$, and alternative and Malcev algebras is infinite \cite{Shestakov77}. The base rank of a variety of algebras with zero multiplication is $0$.

\begin{lm}\label{l4} If $\mathrm{rb}({\mathfrak M})\leq n$, then  the homomorphism (\ref{f4})  is an isomorphism.
\end{lm}
\Proof Suppose that $u$ belongs to the kernel of $\varphi$. This means that $uy_1=0$ is an identity for $A$. Then $uy_1=0$ is an identity for the variety of algebras $\mathfrak{M}$ since $\mathrm{rb}({\mathfrak M})\leq n$. This means $u=0$. $\Box$

\begin{lm}\label{l5} Let ${\mathfrak M}$ be a Nielson-Shreier variety of algebras and
 $D$ be an arbitrary element of $\mathscr{L}(A)$. Then the Jacobian matrix $J(D)$ of $D$ is nilpotent if and only if $R_D$ is a nilpotent element of $M(\mathscr{L}(A))$.
\end{lm}
\Proof Using the proof of Corollary \ref{c1}, it is sufficient to prove that the equality $J(G)^mF=0$ for all $F$ implies that $J(G)^m=0$. This is certainly true if the homomorphism (\ref{f4})  is an isomorphism. Consequently, it is true if $\mathrm{rb}({\mathfrak M})\leq n$ by Lemma \ref{l4}.

Assume that $\mathrm{rb}({\mathfrak M})> n$. The base rank of a (nontrivial) Nielson-Schreier variety of algebras is $2$ if a free algebra in one variable is one dimensional and $1$ otherwise \cite{Um11}. Consequently, we have to consider only the case $\mathrm{rb}({\mathfrak M})=2> n=1$. In this case $A=kx_1$ and $\Omega(\a x_1)=\a$. If $G=(\a x_1)$, then $J(G)=[\a]$.
If $F=(x_1)$, then $J(G)^mF=(\a^mx_1)=0$ implies $\a=0$. $\Box$

The same proof can be applied for the varieties of associative, Novikov, and dual Leibniz algebras.

Thus, the study of left and right nilpotent elements is related to the study of well known combinatorial problems on free algebras. Some examples and counterexamples are given in \cite{UU2014-3} for polynomial algebras. In the case of Nielsen-Schreier varieties we can expect more valuable results. Recall that all automorphisms of finitely generated free algebras are tame in this case \cite{Lewin}.

\begin{prob}\label{pr1} Let ${\mathfrak M}$ be an arbitrary Nielsen-Schreier variety of $m$-ary algebras over $k$ and let $A=k_{\mathfrak M}\langle x_1,x_2,\ldots,x_n\rangle$ be the free algebra of ${\mathfrak M}$ in the variables $x_1,x_2,\ldots,x_n$.

(i) Is every left or right nilpotent element of $\mathscr{L}(A)$ triangulable?

(ii) Is the Jacobian matrix of every right nilpotent derivation nilpotent?

(iii) Is left or right nilpotent element of $\mathscr{L}(A)$ algorithmically recognizable?

(iv) Is every derivation with nilpotent Jacobian matrix triangulable?
\end{prob}

Direct calculations give that the Jacobian matrix of every right nilpotent element of the left-symmetric Witt algebra $\mathscr{L}_2$ is nilpotent.

Problem \ref{pr1}(i) for left nilpotent elements is a well known open problem for free Lie algebras:
Is every locally nilpotent derivation of a finitely generated free Lie algebra triangulable?
It is well known that if $D$ is a derivation of $A$ then $\varphi^{-1}D\varphi$ is a derivation of $A$. But conjugation by an automorphism is not an automorphism of $\mathscr{L}(A)$.
The nilpotency of a Jacobian matrix is also not invariant under conjugation. Using the tameness of automorphisms of $A$ for Nielsen-Schreier varieties, it is possible to show that Problem \ref{pr1}(iv) has a positive solution for homogeneous derivations.

\section{Novikov algebras of derivations}

\hspace*{\parindent}

In this section we examine the left-symmetric algebras of derivations of some free algebras in one variable.

{\bf Example 4.} First we recall the standard basis  of the left-symmetric Witt algebra $\mathscr{L}_1$. Let $k[x]$ be a polynomial algebra in one variable $x$. Put $e_i=x^{i+1}\partial_x$ for all $i\geq -1$.  Then
\bes
\mathscr{L}_1=ke_{-1}\oplus ke_0\oplus ke_1\oplus\ldots \oplus ke_n\oplus\ldots
\ees
and
\bes
e_se_t=(t+1)e_{s+t}, \ \ \ \ s,t\geq -1.
\ees

The algebra $\mathscr{L}_1$ is a {\em Novikov} algebra, i.e., it satisfies also the identity
\bee\label{f5}
(xy)z=(xz)y.
\eee
The positive part ${\mathscr{L}_1}_+$ of $\mathscr{L}_1$ is generated by $e_1$. 

{\bf Example 5.} Recall that an algebra
with a bilinear operation $[-,-]$ is called (right) {\em Leibniz} if it satisfies the {\em Leibniz identity}
\bes
[x,[y,z]]=[[x,y],z]-[[x,z],y].
\ees
Let $B$ be the free Leibniz algebra in one variable $x$.
It is well known that the elements of the form
\bes
x, x^{[2]}=[x,x],\ldots, x^{[n]}=[x^{[n-1]},x], \ldots
\ees
form a linear basis for $B$. Put $f_i=x^{[i+1]}\partial_x$ for all $i\geq 0$. Then
\bes
\mathscr{L}(B)=kf_0\oplus kf_1\oplus\ldots \oplus kf_n\oplus\ldots
\ees
and
\bes
f_0f_i=(i+1)f_i, \ f_sf_t=f_{s+t}, \ \ i,t\geq 0, s\geq 1.
\ees
The positive subalgebra $\mathscr{L}(B)_+$ is generated by $f_1$ and is an associative and commutative algebra. Despite this, $\mathscr{L}(B)$ is not a Novikov algebra.

{\bf Example 6.} An algebra is called {\em
dual Leibniz} if it satisfies the identity
\bes
(xy)z=x(zy+yz).
\ees
Let $C$ be the free dual Leibniz algebra in one variable $x$.
It is well known \cite{Loday95} that the elements of the form
\bes
x, x^2,\ldots, x^n=xx^{n-1}, \ldots
\ees
form a linear basis for $C$ and \cite{Dzhuma05} that
\bee\label{f6}
x^ix^j=
\left(\begin{array}{cc}
i+j-1\\
j\\
\end{array}\right) x^{i+j}
\eee
for all $i,j\geq 1$.

Put $g_i=(i+1)! x^{i+1}\partial_x$ for all $i\geq 0$. Then
\bes
\mathscr{L}(C)=kg_0\oplus kg_1\oplus\ldots \oplus kg_n\oplus\ldots
\ees
and
\bes
g_sg_t=(t+1)g_{s+t}, \ \ \  s,t\geq 0.
\ees
Consequently, $\mathscr{L}(C)$ is a subalgebra of $\mathscr{L}_1$

\begin{pr}\label{p2} Let $A$ be a free (binary) algebra in one variable $x$ with an identity element. Then $\mathscr{L}(A)$ is a Novikov algebra if and only if $A\simeq k[x]$ and $\mathscr{L}(A)\simeq \mathscr{L}_1$.
\end{pr}
\Proof Assume that $\mathscr{L}(A)$ is a Novikov algebra. It is sufficient to prove that $A$ is associative since $A$ is generated by only one element. In the proof of Proposition \ref{p1}, we noticed that the only free associative and commutative algebras with an identity are polynomial algebras. Consequently, $A\simeq k[x]$.

Recall that an algebra is called {\em power-associative} if its every subalgebra generated by one element is associative. It is well known \cite{Albert48} that a variety of algebras over a field of characteristic zero  is power-associative if and only if it satisfies the identities
\bes
x^2x=x^3, x^2x^2=x^4.
\ees

For any $a\in A$ denote by $D_a$ the derivation of $A$ defined by $D_a(x)=a$. By (\ref{f5}), we get
\bes
(D_1D_{x^3})D_{x^2}=(D_1D_{x^2})D_{x^3}.
\ees
Notice that
\bes
(D_1D_{x^3})D_{x^2}=3D_{x^2}D_{x^2}=3D_{x^2x+x^3}, \\
(D_1D_{x^2})D_{x^3}=2D_xD_{x^3}=6D_{x^3}.
\ees
Consequently, $x^2x=x^3$. Using
\bes
(D_1D_{x^4})D_{x^2}=(D_1D_{x^2})D_{x^4},
\ees
we similarly get that $x^3x=x^4$. Consider one more corollary
\bes
(D_xD_{x^2})D_{x^3}=(D_xD_{x^3})D_{x^2}.
\ees
of (\ref{f5}).
Using the preceding equalities we obtain that
\bes
(D_xD_{x^2})D_{x^3}=2D_{x^2}D_{x^3} =2D_{x^2x^2+x(x^2x)+x^4}=2D_{x^2x^2+2x^4},\\
(D_xD_{x^3})D_{x^2}=3D_{x^3}D_{x^2}=3D_{x^3x+x^4}=6D_{x^4}.
\ees
Consequently,
$x^2x^2=x^4$.
$\Box$

If $A$ is a commutative free (binary) algebra in one variable $x$ without an identity element and $\mathscr{L}(A)$ is Novikov algebra, then $A$ is a nilpotent associative and commutative algebra. This can be proved by the same discussions as in the proof of Proposition \ref{p2}.

\section{Derivations of symmetric algebras}

\hspace*{\parindent}

An algebra $A$ over $k$ with one $m$-ary multilinear operation $\langle \cdot,\cdot,\ldots, \cdot\rangle$ is called {\em  symmetric} if
\bes
\langle x_1,x_2,\ldots,x_m\rangle=\langle x_{\sigma(1)},x_{\sigma(2)},\ldots,x_{\sigma(m)}\rangle
\ees
for all $x_1,x_2,\ldots,x_m\in A$ and for any $\sigma\in S_m$, where $S_m$ is the symmetric group on $m$ symbols.
If $m=2$, then a symmetric algebra becomes a commutative (non-associative) algebra.

The main source  of $m$-ary symmetric algebras are polarization algebras \cite{UU2014-2}.
The study of polarization algebras and free $m$-ary symmetric Engel algebras in one variable is motivated by an approach to the Jacobian Conjecture by L.M. Dru\.zkowski and K. Rusek \cite{DR85}, G. Gorni and G. Zampieri \cite{GZ96}, and A.V. Yagzhev \cite{Yagzhev00-1} (see, also \cite{BBRY,UU2014-2}). These results also can be expressed in the language of left-symmetric Witt algebras \cite{UU2014-3}. Recall that the study of commutative (binary symmetric) Engel algebras is related to the study of the well known Alberts problem \cite{Albert}.

\begin{pr}\label{p3} Let $A$ be the free algebra in one variable $x$ of a variety  of symmetric $m$-ary algebras over a field $k$ of characteristic $0$. Then $\mathscr{L}(A)_+$ is generated by $D=\langle x,x,\ldots,x\rangle\partial_x$.
\end{pr}
\Proof Following \cite{KBKA}, we give a definition of all nonassociative words in one letter $x$ with respect to $m$-ary operation $\langle\cdot,\cdot,\ldots,\cdot\rangle$:

(i) $x$ is a unique nonassociative word of length $1$;

(ii) If $w_1,w_2,\ldots,w_m$ are nonassociative words of length $k_1,k_2,\ldots,k_m\geq 1$, respectively, then
$\langle w_1,w_2,\ldots,w_m\rangle$ is a nonassociative word of length $k_1+k_2+\ldots+k_m$.

It is well known \cite{KBKA} that every nonassociative word $w$ of length $l(w)>1$ has a unique decomposition $w=\langle w_1,w_2,\ldots,w_m\rangle$.

Let $u,v$ be  two arbitrary nonassociative words. Put $u<v$ if $l(u)<l(v)$. If
$l(u)=l(v)>1$, then put $u<v$ if $u=\langle u_1,u_2,\ldots,u_m\rangle$, $v=\langle v_1,v_2,\ldots,v_m\rangle$, and $u_1=v_1,u_2=v_2,\ldots,u_{i-1}=v_{i-1}$ and  $u_i<v_i$ for some $1\leq i\leq n$.

A nonassociative word $w$ is called {\em reduced} if $w=\Lambda$  or $w=x$ or $w=\langle w_1,w_2,\ldots,w_m\rangle$ and $w_1\geq w_2\geq\ldots\geq w_m$ are reduced words. Obviously, every element of $B$ is a linear combination of reduced words \cite{Shirshov}.

Let $S=k x\partial_x+\langle D\rangle$ be the subalgebra of $\mathscr{L}(B)$ generated by $D$ and the right identity $x\partial_x$ . We prove that $w\partial_x\in S$ for every reduced word $w$ of length $l(w)\geq 1$. Every reduced word $w$ of length $l(w)> 1$ can be written uniquely as $w=\langle w_1,w_2,\ldots,w_m\rangle$, where $w_1\geq w_2\geq\ldots\geq w_m$ are reduced words of length $\geq 1$. We put $\rho(w)=i$  if $l(w_i)>1$ and $w_{i+1}=\ldots=w_n=x$. Notice that $\rho(w)$ is not defined only if $w_i=x$ for all $i$. In this case we have  $w\partial_x=D\in S$.

Let $w$ be a reduced word with minimal $(l(w),\rho(w))$ such that $w\partial_x\notin S$ (pairs $(l(w),\rho(w))$ are compared lexicographically). Put $\rho(w)=i$. Then $w=\langle w_1,w_2,\ldots,w_m\rangle$, where $w_1\geq w_2\geq\ldots\geq w_m$ are reduced words of length $\geq 1$, $l(w_i)>1$, and $w_{i+1}=\ldots=w_n=x$. Put $u=\langle w_1,w_2,\ldots,w_{i-1},x,\ldots,x\rangle$. Note that $l(u)<l(w)$ and $w_i\partial_x, u\partial_x\in S$ by the choice of $w$. Consequently, $w_i\partial_x u\partial_x=(n-i+1)w\partial_x+t\partial_x\in S$, where $t$ is a linear combination of nonassociative words $u$ such that $l(u)=l(w)$ and $\rho(u)<\rho(w)$. Again $t\partial_x\in S$ by the choice of $w$. Consequently, $w\partial_x\in S$, i.e., a contradiction. $\Box$

We say that an $m$-ary symmetric algebra $A$ is
an {\em Engel} algebra if there exists a positive integer $q$ such that
\bee\label{f7}
M(a,a,\ldots,a)^q b=0
\eee
for all $a,b\in A$. If $m=2$, then this is the standard definition of a commutative Engel algebra.
We also say that a variety of algebras satisfies (\ref{f7}) if every algebra of this variety satisfies (\ref{f7}).

Let $A$ be the free algebra in one variable $x$ of a variety  of symmetric $m$-ary algebras over a field $k$ of characteristic $0$. Consider the endomorphism
\bes
\varphi :A\longrightarrow A, \ \ \ \varphi(x)=x-\langle x,x,\ldots,x\rangle. 
\ees
If $\varphi$ is an automorphism, then the Jacobian matrix $J(D)$ of
 $D=\langle x,x,\ldots,x\rangle\partial_x$ is nilpotent \cite{UU2014-2}. Notice that 
 \bes
 J(D)=[mU(x,\ldots,x)].
 \ees
 Consequently, $J(D)$ is nilpotent if and only if $\mathfrak{M}$ is an Engel variety of algebras. In order to solve the Jacobian Conjecture, it is sufficient to prove that $\varphi$ is an automorphism when $\mathfrak{M}$ is an arbitrary variety of algebras generated by a finite dimensional ternary symmetric Engel algebra \cite{Yagzhev00-1} (see also \cite{UU2014-2,BBRY}). 
 
 By Corollary \ref{c1}, $R_D$ ia a nilpotent element of $\mathscr{L}(A)$ and $D$ is a right nilpotent element of $D$. As in \cite{UU2014-3}, the coefficients of the formal inverse to $\varphi$ can be expressed as Lie polynomials on right powers $D^{[i]}$ of $D$. Problem \ref{pr1}(ii) asks: Does the right nilpotency of $D$ imply the nilpotency of $R_D$. Here we positively solve some elementary cases.

\begin{lm}\label{l6}
Let $A$ be the free algebra in one variable $x$ of a variety  of symmetric $m$-ary algebras over a field $k$ of characteristic $0$ and let $D=\langle x,x,\ldots,x\rangle\partial_x$. Then the following statements are true:

$(i)$ If $D^2=0$, then $R_D$ is nilpotent;

$(ii)$ If $m=2$ and $D^{[3]}=0$, then $R_D$ is nilpotent.
\end{lm}
\Proof $(i)$ Assume that $D^2=0$. This implies that
\bee\label{f8}
\langle x,x,\ldots,x,\langle x,x,\ldots,x\rangle\rangle=0.
\eee

Put $w=\langle x,x,\ldots,x\rangle$. Then $D=w\partial_x$ and (\ref{f8}) means that $D(w)=0$. We prove that $A$ is the linear span of $x$ and $w$. For this, it is sufficient to prove that
\bee\label{f9}
\langle x,x,\ldots,x,\underbrace{w,\ldots,w}_i\rangle=0
\eee
for all $1\leq i\leq m$. If $i=1$, then (\ref{f9}) coincides with (\ref{f8}). Suppose that
(\ref{f9}) is true for some $1\leq i<m$. Applying $D$ to  (\ref{f9}), we get
\bes
(m-i) \langle x,x,\ldots,x,\underbrace{w,\ldots,w}_{i+1}\rangle=0
\ees
since $D(w)=0$.

Thus, $A$ is nilpotent and $\varphi$ is an automorphism of $A$. This means that $J(D)$ is nilpotent.

$(ii)$. In this case $A$ is a binary commutative algebra and satisfies the identity
\bee\label{f10}
x(x(xx))=0.
\eee
The product $a_1(a_2(\ldots(a_{m-1}a_m)\ldots))$ will be denoted by $a_1a_2\ldots a_m$.
A partial linearization of (\ref{f10}) gives
\bee\label{f11}
yxxx+xyxx+2xxxy=0.
\eee
This identity can be written in the form
\bee\label{f12}
L_{xx^2}+L_xL_{x^2}+2L_xL_xL_x=0
\eee
Applying $D$ to this identity, we also get
\bee\label{f13}
L_{x^2x^2}+L_{x^2}L_{x^2}+2L_xL_{xx^2}+2L_{x^2}L_xL_x+2L_xL_{x^2}L_x+2L_xL_xL_{x^2}=0.
\eee
Direct calculations, using the identities (\ref{f10})--(\ref{f13}), give that
\bes
x^2xx^2=-xx^2x^2, x^2x^2x^2=xxx^2x^2, (xx^2)xx^2=xxx^2x^2,\\
 xxxx^2x^2=0, x^2xx^2x^2=0, (xx^2)x^2x^2=0,
\ees
and $A$ is a linear span of
\bes
x, x^2, xx^2, x^2x^2,
xx^2x^2,
xxx^2x^2.
\ees
Consequently, $A$ is nilpotent and $\varphi$ is an automorphism, i.e., $J(D)$ is nilpotent.
$\Box$

\end{document}